\documentstyle[12pt]{amsart}

\newtheorem{theo+}              {Theorem}           [section]
\newtheorem{prop+}  [theo+]     {Proposition}
\newtheorem{coro+}  [theo+]     {Corollary}
\newtheorem{lemm+}  [theo+]     {Lemma}
\newtheorem{exam+}  [theo+]     {Example}
\newtheorem{rema+}  [theo+]     {Remark}
\newtheorem{defi+}  [theo+]     {Definition}

\newenvironment{theorem}{\begin{theo+}}{\end{theo+}}
\newenvironment{proposition}{\begin{prop+}}{\end{prop+}}

\newenvironment{lemma}{\begin{lemm+}}{\end{lemm+}}
\newenvironment{example}{\begin{exam+}}{\end{exam+}}
\newenvironment{remark}{\begin{rema+}}{\end{rema+}}
\newenvironment{definition}{\begin{defi+}}{\end{defi+}}

\begin{document}
\baselineskip  18pt

\title [Complete lifts of harmonic maps and morphisms...]{Complete lifts of
harmonic maps and morphisms between Euclidean spaces}
\bigskip
\thanks{*The author was supported by the Chinese Education Commission \\
as ``Visiting scholar program (94A04)''.}
\maketitle

\begin{center}
{\bf Ye-lin Ou*}\\
{Department of Mathematics, Guangxi University for Nationalities, Nanning
530006, P.R.China}\\
Current : Department of Pure Mathematics, School of Mathematics, University of
Leeds, Leeds LS2 9JT, U.K. (until Jan.1996). \\
E-mail: pmtyo@@amsta.leeds.ac.uk\\
\end{center}
\bigskip

\begin{center}
{\bf Abstract}
\end{center}
\bigskip
We introduce the complete lifts of maps between (real and complex) Euclidean
spaces and study their properties concerning holomorphicity, harmonicity and
horizontal weakly conformality. As applications, we are able to use this
concept to characterize holomorphic maps  $\phi:{\Bbb C}^{m}\supset
U\longrightarrow {\Bbb  C}^{n}$ (Proposition 2.3) and  to construct many new
examples of harmonic morphisms (Theorem 3.3). Finally we show that the
complete lift of the quaternion product followed by the complex product is a
simple and explicit example of a harmonic morphism which does not arise (see
Definition 4.8 in \cite{BaiWoo95}) from any K{\"a}hler structure.
\bigskip


\newpage

\section{Introduction}

Let$({M}^{m},g)$ and $({N}^{n},h)$ be two Riemannian manifolds. A map
$\phi:({M}^{m},g)\longrightarrow ({N}^{n},h)$ is called a harmonic map if the
divergence of its differential vanishes. Such maps are the extrema  of the
energy functional $ \frac{1}{2}\int_{D} {\left|d\phi\right|}^2$ over compact
domain $D\subset M$. For a detailed account on harmonic maps we refer to
\cite{EelLem78},\cite{EelLem83},\cite{EelLem88} and the references therein.
Harmonic morphisms are a special subclass of harmonic maps which preserve
solutions of Laplace's equation in the sense that for any harmonic function
$f:U \longrightarrow{\Bbb R}$, defined on an open subset U of $N$ with
${\phi}^{-1}(U) $ non-empty, $f\circ\phi :{\phi}^{-1}(U) \longrightarrow {\Bbb
R} $ is a harmonic function. In other words, $\phi$ pulls back germs of
harmonic functions on $N$ to germs of harmonic functions on $M$. In the theory
of stochastic processes, harmonic morphisms $\phi:(M,g)\longrightarrow (N,h)$
are found t!
 o be Brownian path preserving mappings meaning that they map Brownian motions
on $M$ to Brownian motions on $N $(\cite{BerCamDav79},\cite{Levy}). The
following nice characterization of a harmonic morphism is due to Fuglede and
Ishihara independently.
\begin{theorem}
(\cite{Fug78},\cite{Ish79A})A map $\phi$ is a harmonic morphism if and only if
it is both a harmonic and a horizontal weakly conformal map.
\end{theorem}
For a map  $\phi:{\Bbb  R}^{m}\supset U\longrightarrow{\Bbb R}^{n}$ between
Euclidean spaces, with $\phi(x)=({\phi}^{1}(x),...,{\phi}^{n}(x))$, the
harmonicity and weakly conformality are just equivalent to the following
conditions respectively:
\begin{align}
&\sum_{i=1}^{m} \frac{\partial^{2}{\phi}^{k}}{\partial{x}_{i}^{2}} = 0 \\
&\sum_{i=1}^{m} \frac{\partial{\phi}^{k}}{\partial{x}_{i}}
\frac{\partial{\phi}^{l}}{\partial{x}_{i}} = {\lambda}^{2}(x) {\delta}_{kl}
\end{align}
For $k,l =1,2,...,n.$

 It is well-known that (i) if $ \mbox{dim}\;M < \mbox{dim}\;N$ then any
harmonic morphism is a constant map, (ii) if $\mbox{dim}\;M = \mbox{dim}\;N
=2$, the harmonic morphisms are precisely the weakly conformal maps, (iii) if
$\underline{\mbox{dim} \;N = 1}$ they are precisely the harmonic maps, (iv) if
$\mbox{dim}\;M = \mbox{dim}\;N \geq 3$ they are precisely the homothetic
maps(\cite{Fug78}, \cite{Ish79A}).
Though the composition of harmonic maps does not generally turn out to be a
harmonic map , we have
\begin{lemma}
If ${\phi}_{1} :(M,g)\longrightarrow(\overline{M},\overline{g})$ and $
{\phi}_{2}: (\overline{M},\overline{g}) \longrightarrow (N,h)$ are harmonic
morphisms, then so is their composition
${\phi}_{2}\circ{\phi}_{1}:(M,g)\longrightarrow(N,h).$
\end{lemma}

\begin{lemma}{(\cite{Gud94A})}
Let $\pi :(\overline{M},\overline{g})\longrightarrow(M,g)$ be a non-constant
harmonic morphism, $\phi : (M,g)\longrightarrow (N,h)$ be a map, and $\Phi
=\phi\circ\pi :(\overline{M},\overline{g})\longrightarrow (N,h)$ their
composition, then $\Phi$ is a harmonic morphism if and only if $\phi$ is a
harmonic morphism.
\end{lemma}
Since the equations determining a harmonic morphism form an overdetermined
system of partial differential equations one should not expect to find too many
examples of harmonic morphisms, however, one can easily check that all the maps
in the following example are in this class.
\begin{example}
(i)  $\phi:{\Bbb C}^{2}\longrightarrow {\Bbb C}, \phi(z,w) = zw$ and $\phi(z,w)
= z {\overline{w}}$;\\
(ii) $\phi :{\Bbb C}^{2}\longrightarrow{\Bbb R}\times{\Bbb C}$, with $
\phi(z,w) = ({| z|}^{2}-{| w|}^{2},2zw)$;\\
(iii) The quaternion product \\$q :{\Bbb C}^{2}\times{\Bbb
C}^{2}\longrightarrow{\Bbb C}^{2}, q({z}_{1},{z}_{2},{z}_{3},{z}_{4}) =(
{z}_{1}{z}_{3}-{z}_{2}{\overline{z}}_{4},{z}_{1}{z}_{4}+{z}_{2}{\overline{z}}_{3})$.\\
(iv) The hyperbolic analogue of stereographic projection \\
$\phi:{\Bbb R}^{3}\backslash  \left\{(0,0,{x}_{3})\mid{x}_{3}\geq
0\right\}\longrightarrow{\Bbb R}^{2}$ with $\phi({x}_{1},{x}_{2},{x}_{3})=
\left(\frac{{x}_{1}}{r-{x}_{3}},\frac{{x}_{2}}{r-{x}_{3}}\right)$, where
${r}^{2}={x}_{1}^{2}+{x}_{2}^{2}+{x}_{3}^{2}$.\\
(v) The orthogonal projection $\phi:{\Bbb R}^{m}\longrightarrow{\Bbb R}^{n}$
with $\phi({x}_{1},\ldots,{x}_{m}) =({x}_{1},\ldots,{x}_{n})$, and the radial
projection $\phi : {\Bbb R}^{m}\backslash \{0\} \longrightarrow{S}^{m-1}
(m=1,2,\ldots) $ with $\phi(x) = x/\left|x\right|$.\\
(vi) The natural projection $\pi: (TM,\overline{g})\longrightarrow(M,g)$, where
$\overline{g}$ is the Sasaki metric on the tangent bundle of the Riemannian
manifold $(M,g)$.\\
(vii) the Hopf maps ${S}^{3}\longrightarrow{S}^{2}$ and
${S}^{2n+1}\longrightarrow{\Bbb C}{P}^{n} (n = 1,2,\ldots).$\\
(viii) Any holomorphic maps from a K{\"a}hler manifold to a Riemann surface
(\cite{BaiEel81}), in particular, any holomorphic function $\phi :
U\longrightarrow{\Bbb C}$ on an open subset $U$ of ${\Bbb C}^{m}$ is a harmonic
morphism.
\end{example}
In recent years, much work has been done in classifying and constructing
harmonic morphisms from certain model spaces to other model spaces (see e.g.
\cite{BaiWoo88}, \cite{BaiWoo91}, \cite{BaiWoo92A} ,\cite{BaiWoo95},
\cite{Gud92}, \cite{Gud94A} , \cite{Gud94},  \cite{GudSig93}, \cite{Woo86A},
\cite{Woo92}). More recently Baird and Wood have found (see \cite{Woo92}
,\cite{BaiWoo95}) some strong links between Hermitian structures and harmonic
morphisms from open subsets of ${\Bbb R}^{2n}$ to ${\Bbb C}$ or a Riemann
surface. They have constructed many interesting locally and globally defined
harmonic morphisms. Some of them are holomorphic with respect to non-K{\"a}hler
structures and some are not holomorphic with respect to any K{\"a}hler
structure. In this work we use the complete lift of the quaternion product and
the composition with the complex product to give a harmonic morphism
$\Phi:{\Bbb R}^{16}\supset U\longrightarrow {\Bbb  C}$ which is not holomorphic
with respe!
 ct to any K{\"a}hler structure.

\section{Complete Lifts and Their Properties}

\begin{definition}
Let $\phi:{\Bbb R}^{m}\supset U\longrightarrow {\Bbb  R}^{n},
\phi(x)=({\phi}^{1}(x),...,{\phi}^{n}(x))$, be a map from an open connected
subset of   ${\Bbb R}^{m} $ into $ {\Bbb  R}^{n}$. The (real) complete lift of
$\phi$ is a map $\Phi:{\Bbb R}^{2m}\supset {U\times{\Bbb R}^{m}}\longrightarrow
{\Bbb  R}^{n}$, given by
\begin{equation}
\Phi({x}_{1},\ldots,{x}_{m};{y}_{1},\ldots,{y}_{m})
=\left(\frac{\partial{\phi}^{i}}{\partial{x}_{j}}(x)\right)
\left(\begin{array}{c}
y_{1}\\\vdots\\y_{m}
\end{array}\right)
\end{equation}
where $\left(\frac{\partial{\phi}^{i}}{\partial{x}_{j}}(x)\right)$ denotes the
Jacobian matrix of $\phi$ at $x$.
\end{definition}

\begin{remark}
(i) The complete lift of a map is a partial linear map in the sense that it
depends linearly on half of its variables.\\(ii) Let  $\phi:{\Bbb C}^{m}\supset
U\longrightarrow {\Bbb  C}^{n}$ be a ${C}^{\infty}$ from an open connected
subset of ${\Bbb C}^{m}$ into ${\Bbb  C}^{n}$ ,then the (complex) complete lift
of $\phi$ can be defined similarly by using the complex Jacobian matrix.
\end{remark}
It is well-known that ${\Bbb C}^{m}$ can be identified with ${\Bbb R}^{2m}$ by
identifying $({z}_{1},\ldots,{z}_{m})$ with
$({x}_{1},{y}_{1},\ldots,{x}_{m},{y}_{m})$, where ${z}_{k} = {x}_{k} +
i{y}_{k}.$  Thus any map  $\phi:{\Bbb C}^{m}\supset U\longrightarrow {\Bbb
C}^{n}$  can be identified with a map (called {\bf the real identification of}
$\phi$)  ${\phi}_{r}:{\Bbb R}^{2m}\supset W\longrightarrow {\Bbb  R}^{2n}$,
given by
\begin{equation}
 {\phi}_{r}({x}_{1},{y}_{1},\ldots,{x}_{m},{y}_{m}) = (
{u}^{1}(x,y),{v}^{1}(x,y),\ldots,{u}^{n}(x,y),{v}^{n}(x,y))\notag
\end{equation}
where ${u}^{k}(x,y)$and ${v}^{k}(x,y)$  are the real and imaginary parts of the
 $ kth $ component function of $\phi \,$ , i.e  $ \,{\phi}^{k} (x,y) =
{u}^{k}(x,y) +i{v}^{k}(x,y)$.\\

It should be noted that the operation of the complete lift and the above
identification does not commute in general. In fact, the following is easily
established:
\begin{proposition}
Let $\phi : U\longrightarrow{\Bbb C}^{n} $ be a ${C}^{\infty}$ map from an open
connected subset $U$ of ${\Bbb C}^{m}$, then  $\phi$ is holomorphic if and only
if the (complex ) complete lift $\Phi$ of $\phi$ is identical with the (real)
complete lift ${\Phi}_{r}$ of its real identification ${\phi}_{r}$, viewed as a
map between complex Euclidean spaces.
\end{proposition}

\begin{example}
The quaternion product $q :{\Bbb C}^{2}\times{\Bbb C}^{2}\longrightarrow{\Bbb
C}^{2}, q({z}_{1},{z}_{2},{z}_{3},{z}_{4}) =(
{z}_{1}{z}_{3}-{z}_{2}{\overline{z}}_{4},{z}_{1}{z}_{4}+{z}_{2}{\overline{z}}_{3})$ is not holomorphic. The  complex complete lift of $ Q $ is given by
\begin{equation}
Q ({z}_{1},{z}_{2},{z}_{3},{z}_{4},{w}_{1},{w}_{2},{w}_{3},{w}_{4}) =
({z}_{3}{w}_{1} - {\overline{z}}_{4}{w}_{2} + {z}_{1}{w}_{3} , {z}_{4}{w}_{1} +
{\overline{z}}_{3}{w}_{2} + {z}_{1}{w}_{4})\notag
\end{equation}
whilst the real complete lift of $q_{r}$, viewed as a complex map, is given by
\begin{align}
Q_{r} &({z}_{1},{z}_{2},{z}_{3},{z}_{4},{w}_{1},{w}_{2},{w}_{3},{w}_{4}) =
\notag\\
      &({z}_{3}{w}_{1} - {\overline{z}}_{4}{w}_{2} + {z}_{1}{w}_{3} -
{z}_{2}{\overline{w}}_{4}, {z}_{4}{w}_{1} + {\overline{z}}_{3}{w}_{2} +
{z}_{2}{\overline{w}}_{3} + {z}_{1}{w}_{4})\notag
\end{align}
which is obviously different from  $ Q$ .
\end{example}

\begin{theorem}
The complete lift of any holomorphic map  \\
$\phi:{\Bbb C}^{m}\supset U\longrightarrow {\Bbb  C}^{n}$ is holomorphic .
\end{theorem}

\begin{proof}
It suffices to check that the component function  ${\Phi}^{l}(z,w) =
\sum_{k=1}^{m}\frac{\partial{\phi}^{l}}{\partial{z}_{k}}{w}_{k} ( l =
1,2,\ldots,n) $ is holomorphic with respect to  variables
${z}_{1},\ldots,{z}_{m},{w}_{1},\ldots,{w}_{m}$. This amounts to check that the
following equations
\begin{align}
&\frac{\partial}{\partial{{\overline{z}}_{k}}}(\sum_{j=1}^{m}\frac{\partial{\phi}^{l}}{\partial{z}_{j}}{w}_{j}) = 0 \notag\\
&\frac{\partial}{\partial{{\overline{w}}_{k}}}(\sum_{j=1}^{m}\frac{\partial{\phi}^{l}}{\partial{z}_{j}}{w}_{j}) = 0 \notag\
\end{align}
hold for $ l =1,2,\ldots,n,\; k = 1, 2,\ldots, m$ which is trivial and is
omitted.
\end{proof}

\begin{theorem}
The complete lift of any harmonic map  $\phi:{\Bbb R}^{m}\supset
U\longrightarrow {\Bbb  R}^{n}$ is harmonic.
\end{theorem}

\begin{proof}
Let $\phi:{\Bbb R}^{m}\supset U\longrightarrow {\Bbb  R}^{n},
\phi(x)=({\phi}^{1}(x),...,{\phi}^{n}(x))$, be a harmonic map, then we have
$\sum_{j=1}^{m} \frac{\partial^{2}{\phi}^{k}}{\partial{x}_{j}^{2}} = 0 $ for $
k = 1,2,\ldots,n.$ We must check that the component function of the complete
lift  ${\Phi}^{k}(x,y) =  \sum_{i=1}^{m}
\frac{\partial{\phi}^{k}}{\partial{x}_{i}}{y}_{i} $,  $ k = 1,2,\ldots,n$ , is
a harmonic function of ${x}_{1},\ldots,{x}_{m},{y}_{1},\ldots,{y}_{m}$. But\\
$\sum_{j=1}^{m}\frac{\partial^{2}}{\partial{x}_{j}^{2}}(\sum_{i=1}^{m}
\frac{\partial{\phi}^{k}}{\partial{x}_{i}}{y}_{i}) +
\sum_{j=1}^{m}\frac{\partial^{2}}{\partial{y}_{j}^{2}}(\sum_{i=1}^{m}
\frac{\partial{\phi}^{k}}{\partial{x}_{i}}{y}_{i})
=\sum_{i=1}^{m}{y}_{i}\frac{\partial}{\partial{x}_{i}}(\sum_{j=1}^{m}
\frac{\partial^{2}{\phi}^{k}}{\partial{x}_{j}^{2}}) = 0$. This ends the proof.
\end{proof}

For the complete  lift of a horizontal conformal map we have

\begin{theorem}
Let  $\phi:{\Bbb R}^{m}\supset U\longrightarrow {\Bbb  R}^{n}$ be a horizontal
weakly conformal map, then the complete lift $\Phi$  of $\phi$ is  horizantal
weakly conformal if and only if the following conditions hold for $\alpha,\beta
= 1, 2,\ldots, n.$:
\begin{align}
&(\em{hess} \phi^\alpha)^{2} = (\mbox{hess} \phi^\beta)^{2}\\
&(\em{hess}\phi^\alpha)(\mbox{hess}\phi^\beta) = -
(\mbox{hess}\phi^\beta)(\mbox{hess}\phi^\alpha)
\end{align}
 where $(\em{hess}\phi^\alpha)$ denotes the Hessian matrix of the component
function $\phi^\alpha$.
\end{theorem}

\begin{proof}
At each point $(x,y)$, the Jacobian matrix of $\Phi$ is given by
$$
J(\Phi) =
\left( \begin{array}{cccccc}
\partial_{1} \nabla\phi^{1} \bullet y & \ldots& \partial_{m}
\nabla\phi^{1}\bullet y & \frac{\partial{\phi}^{1}}{\partial{x}_{1}}& \ldots
&\frac{\partial{\phi}^{1}}{\partial{x}_{m}}\\

\partial_{1} \nabla\phi^{2} \bullet y & \ldots& \partial_{m}
\nabla\phi^{2}\bullet y & \frac{\partial{\phi}^{2}}{\partial{x}_{1}} & \ldots
&\frac{\partial{\phi}^{2}}{\partial{x}_{m}}\\
\vdots & \vdots & \vdots & \vdots& \vdots & \vdots \\

\partial_{1} \nabla\phi^{n} \bullet y & \ldots & \partial_{m}
\nabla\phi^{n}\bullet y & \frac{\partial{\phi}^{n}}{\partial{x}_{1}} & \ldots &
\frac{\partial{\phi}^{n}}{\partial{x}_{m}}
\end{array} \right )
$$
Where $\bullet$ denotes the inner product in Euclidean space. \\
{}From this together with the horizontal weak conformality of $\phi$ it follows
that $\Phi$ is  horizontally weakly conformal if and only if the following
equations hold
\begin{align}
&\sum_{i=1}^{m}(\partial_{i} \nabla\phi^{\alpha} \bullet y)^{2}\equiv
\sum_{i=1}^{m}(\partial_{i} \nabla\phi^{\beta} \bullet y)^{2}\\
&\sum_{i=1}^{m}(\partial_{i} \nabla\phi^{\alpha} \bullet y)(\partial_{i}
\nabla\phi^{\beta} \bullet y)\equiv 0
\end{align}
Note that Equations (6) and (7) express the identity of quadratics in
$y_{i}'s$, they are equivalent to the following equations
\begin{align}
&\partial_{i} \nabla\phi^{\alpha} \bullet \partial_{i} \nabla\phi^{\alpha}=
\partial_{i} \nabla\phi^{\beta} \bullet \partial_{i} \nabla\phi^{\beta}\\
&\partial_{i} \nabla\phi^{\alpha} \bullet \partial_{j} \nabla\phi^{\alpha}=
\partial_{i} \nabla\phi^{\beta} \bullet \partial_{j} \nabla\phi^{\beta}\\
&\partial_{i} \nabla\phi^{\alpha} \bullet \partial_{i} \nabla\phi^{\beta} = 0
\\
&\partial_{i} \nabla\phi^{\alpha} \bullet \partial_{j} \nabla\phi^{\beta} = -
\partial_{j} \nabla\phi^{\alpha} \bullet \partial_{i} \nabla\phi^{\beta}
\end{align}
Now equations (8) and (9) are equivalent to the fact that the square of the
Hessian matrix of $\phi^{\alpha}$ is the same for all $\alpha = 1,2,\ldots,n$,
while equations (10) and (11) just say that the product matrix
$(\mbox{hess}\phi^{\alpha})(\mbox{hess}\phi^{\beta})$ is skew symmetric. On the
other hand, since $
(\mbox{hess}\phi^{\alpha})$ and $(\mbox{hess}\phi^{\beta})$ are symmetric it is
easily seen that (10) and (11) are, in fact, equivalent to (5) which ends the
proof.
\end{proof}

\section{Examples and further results}
In this section we give examples of harmonic morphisms whose complete lifts are
harmonic morphisms. We prove that the complete lift of a quadratic harmonic
morphism is always a quadratic harmonic morphism. This gives us a large class
of harmonic morphisms between Euclidean spaces. Two counter examples are also
given, and finally we use the complete lift of the quaternion product and the
composition with the complex product to construct an example of harmonic
morphism $\phi : {\Bbb R}^{16}\longrightarrow {\Bbb C} $ which does not arise
from any K{\"a}hler structure.
\begin{example}
(i)  $\phi:{\Bbb C}^{2}\longrightarrow {\Bbb C},$  $ \phi(z,w) = z
{\overline{w}}$  is a non-holomorphic harmonic morphism. The complex and the
real complete lifts are given respectively by\\
 $\Phi(z_{1},z_{2},w_{1},w_{2}) = \overline{z_{2}} w_{2}$, and
$\Phi_{r}(x_{1},\ldots,x_{4},y_{1},\ldots,y_{4}) = (
x_{3}y_{1}+x_{4}y_{2}+x_{1}y_{3}+x_{2}y_{4} ,
-x_{4}y_{1}+x_{3}y_{2}+x_{2}y_{3}-x_{1}y_{4} )$. It is easily checked that they
are both harmonic morphisms.\\
(ii) $\phi :{\Bbb C}^{2}\longrightarrow{\Bbb R}\times{\Bbb C}
 ,  \phi(z,w) = ({| z|}^{2}-{| w|}^{2},2zw)$ , viewed as  a map between real
Euclidean spaces can be written as\\
$\phi (x_{1},\ldots,x_{4}) = \left( x_{1}^{2}+x_{2}^{2}-x_{3}^{2}-x_{4}^{2} ,
2x_{1}x_{3}-  2x_{2}x_{4} ,  2x_{1}x_{4}- 2x_{2}x_{3}\right)$.\\
Then $\phi$ is a harmonic morphism (see \cite{BerCamDav79}) and its complete
lift is given by\\
$\Phi = ( 2 x_{1}y_{1} + 2 x_{2}y_{2} - 2 x_{3}y_{3} - 2 x_{4}y_{4}  ,  2
x_{3}y_{1} - 2 x_{4}y_{2} + 2 x_{1}y_{3} - 2 x_{2}y_{4}  ,   2 x_{4}y_{1} + 2
x_{3}y_{2} + 2 x_{2}y_{3} + 2 x_{1}y_{4})$.\\
One can easily check that $ \Phi$ is also a harmonic morphism.\\
(iii) The real complete lift of the quaternion product is again a harmonic
morphism while the complex complete lift of the quaternion product is no longer
a harmonic morphism. In fact, the real complete lift of $ q$ is given by\\
\begin{equation}
Q_{r}(x,y)=\left(\begin{array}{cccccccc}
x_{5}&-x_{6}&-x_{7}&-x_{8}&x_{1}&-x_{2}&-x_{3}&-x_{4}\\
x_{6}&x_{5}&x_{8}&-x_{7}&x_{2}&x_{1}&-x_{4}&x_{3}\\
x_{7}&-x_{8}&x_{5}&x_{6}&x_{3}&x_{4}&x_{1}&-x_{2}\\
x_{8}&x_{7}&-x_{6}&x_{5}&x_{4}&-x_{3}&x_{2}&x_{1}
\end{array}\right)\left(\begin{array}{c}
y_{1}\\\vdots\\y_{8}
\end{array}\right)\notag
\end{equation}
It follows from Theorem 2.6 that $Q_{r}(x,y)$ is harmonic.  $Q_{r}(x,y)$ is
also horizontal weakly conformal as one can see this at a glance of its
Jacobian matrix\\
\begin{equation}
J(Q_{r}(x,y)) =  \left(\begin{array}{cc}
J(q_{r}(y)) \mid J(q_{r}(x))\end{array}\right)_{(4\times 16)}.\notag
\end{equation}
On the other hand, the complex complete lift $Q$ of $q$ is given by\\
\begin{equation}
Q(z,w) = ({z}_{3}{w}_{1}-{\overline{z}}_{4}{w}_{2}+{z}_{1}{w}_{3}  ,
{z}_{4}{w}_{1}+ {\overline{z}}_{3}{w}_{2}+{z}_{1}{w}_{4}).\notag
\end{equation}
When viewed as a map between Euclidean spaces its Jacobian matrix can be
calculated as\\
\begin{align}
&J(Q(x,y))= \notag\\
&\left(\begin{array}{cccccccccccccccc}
y_{5}&-y_{6}&0&0&y_{1}&-y_{2}&-y_{3}&-y_{4}&x_{5}&-x_{6}&-x_{7}&-x_{8}&x_{1}&-x_{2}&0&0\\
y_{6}&y_{5}&0&0&y_{2}&y_{1}&-y_{4}&y_{3}&x_{6}&x_{5}&x_{8}&-x_{7}&x_{2}&x_{1}&0&0\\
y_{7}&-y_{8}&0&0&y_{3}&y_{4}&y_{1}&-y_{2}&x_{7}&-x_{8}&x_{5}&x_{6}&0&0&x_{1}&-x_{2}\\
y_{8}&y_{7}&0&0&y_{4}&-y_{3}&y_{2}&y_{1}&x_{8}&x_{7}&-x_{6}&x_{5}&0&0&x_{2}&x_{1}
\end{array}\right)\notag
\end{align}
Therefore, the complex complete lift of $ q $ is not horizontal weakly
conformal though it is  a harmonic map.\\
(iv) The complete lift of the hyperbolic analogue of stereographic projection
in Example 1.4 (iv) is not a harmonic morphism.
\end{example}

Note that the maps in $ (i)-(iii)$ of Example 3.1 are among a large class of
harmonic morphisms whose complete lifts are always harmonic morphisms:
\begin{definition}
A map  $\Phi : {\Bbb R}^{m}\longrightarrow  {\Bbb R}^{n} $ is called a
quadratic map if all of its components are homogeneous polynomials of degree
$2$.
\end{definition}
\begin{theorem}
The complete lift of any quadratic harmonic morphism is again a quadratic
harmonic morphism.
\end{theorem}
\begin{proof}
Let  $\phi : {\Bbb R}^{m}\longrightarrow  {\Bbb R}^{n} $ be a quadratic
harmonic morphism, then by definition we can write
\begin{equation}
\phi (X) = ( X^{t}A_{1}X, \ldots, X^{t}A_{n}X )
\end{equation}
where $ X $ denotes the column vector in ${\Bbb R}^{m}$ , $ X^{t}$ the
transpose of  $ X $   and   $A_{i} (i = 1,...,n) $  is a symmetric matrix of $m
\times m $. We can write the Jacobian matrix of  $\phi (X)$  as
\begin{equation}
J(\phi (X))  = \left(\begin{array}{c}
2 X^{t}A_{1}\\ \vdots\\ 2X^{t}A_{n}\end{array}\right)_{n \times m } \notag
\end{equation}
By definition, the complete lift of $\phi$ can be written as
\begin{align}
\Phi (X,Y) & = \left(\begin{array}{c}
2 X^{t}A_{1}\\ \vdots\\ 2X^{t}A_{n}\end{array}\right) Y \notag\\
& = ( 2X^{t}A_{1}Y, \ldots, 2X^{t}A_{n}Y ) \notag
\end{align}
Now the harmonicity of $\Phi (X,Y)$ follows from that of  $\phi$ by Theorem
2.6. On the other hand, a routine calculation gives

\begin{align}
J(\Phi (X,Y))  &= \left(\begin{array}{c|c}
2Y^{t}A_{1} & 2X^{t}A_{1}\\\vdots & \vdots\\  2Y^{t}A_{n} & 2X^{t}A_{n}
\end{array}\right)_{n \times 2m }\\
&= \left( J(\phi (Y)) | J(\phi (X))\right)_{n \times 2m} \notag
\end{align}
where the second equality holds because of the  fact that $ A_{i} (i =
1,\ldots,n)$ is symmetric.\\
{}From equation (13) one can easily see that $\Phi (X,Y) $ is horizontal weakly
conformal if and only if $\phi$ possesses that property. Thus we have proved
the Theorem.
\end{proof}
\begin{remark}
a) A complete classification of quadratic harmonic morphisms between Euclidean
spaces will appear in \cite{OuWoo95}.\\
b) Though the complete lift of a quadratic harmonic morphism is always a
quadratic harmonic morphism there exist quadratic harmonic morphisms which are
not the complete lift of any map as the following example shows.
\end {remark}

\begin{example}
The quaternion product is not the complete lift of any map. In fact, if there
were a map
\begin{equation}
 \phi(x)=({\phi}^{1}(x),...,{\phi}^{4}(x)) : {\Bbb R}^{4}\longrightarrow  {\Bbb
R}^{4} \notag
\end{equation}
whose complete lift were the quaternion product, then $
{\phi}^{1}(x),...,{\phi}^{4}(x) $ would satisfy the following systems of PDEs
:\\

$$\begin{array}{cccc}
\frac{\partial{\phi}^{1}}{\partial{x}_{1}} = x_{1} &
\frac{\partial{\phi}^{1}}{\partial{x}_{2}} = -x_{2} &
\frac{\partial{\phi}^{1}}{\partial{x}_{3}} = -x_{3} &
\frac{\partial{\phi}^{1}}{\partial{x}_{4}} = -x_{4} \\
\frac{\partial{\phi}^{2}}{\partial{x}_{1}} = -x_{2} &
\frac{\partial{\phi}^{2}}{\partial{x}_{2}} = x_{1} &
\frac{\partial{\phi}^{2}}{\partial{x}_{3}} = -x_{4} &
\frac{\partial{\phi}^{2}}{\partial{x}_{4}} = x_{3}\\
\frac{\partial{\phi}^{3}}{\partial{x}_{1}} = x_{3} &
\frac{\partial{\phi}^{3}}{\partial{x}_{2}} = x_{4} &
\frac{\partial{\phi}^{3}}{\partial{x}_{3}} = x_{1} &
\frac{\partial{\phi}^{}}{\partial{x}_{4}} = -x_{2}\\
\frac{\partial{\phi}^{4}}{\partial{x}_{1}} = x_{3} &
\frac{\partial{\phi}^{4}}{\partial{x}_{2}} = x_{4} &
\frac{\partial{\phi}^{4}}{\partial{x}_{3}} = x_{2} &
\frac{\partial{\phi}^{4}}{\partial{x}_{4}} = x_{1}
\end{array}$$
But this is impossible because, for instance,
$\frac{\partial^{2}\phi^{2}}{\partial x_{2} \partial x_{1}} = -1 \neq
\frac{\partial^{2}\phi ^{2}}{\partial x_{1} \partial x_{2}} = 1 $ .
\end{example}

\begin{example}{Harmonic morphism  $\Phi : {\Bbb R}^{16}\longrightarrow {\Bbb
C} $  not arising from any K{\"a}hler structure.}
We have seen that $Q_{r}(x,y): {\Bbb R}^{16}\longrightarrow {\Bbb C}^{2}$ is a
harmonic morphism, and it is known that  $\phi:{\Bbb C}^{2}\longrightarrow
{\Bbb C}, \phi(z,w) = zw$ is a harmonic morphism. We claim that $\Phi : {\Bbb
R}^{16}\longrightarrow {\Bbb C} $  given by $\Phi = \phi\circ Q_{r}$ is a
harmonic morphism not arising from any K{\"a}hler structure. To see this we
need the following criterion\\

\begin{proposition}{(cf.Baird and Wood \cite{BaiWoo95})}
Let $\Phi : {\Bbb R}^{2m}\supset U\longrightarrow {\Bbb C} $ be a submersive
harmonic morphism. Then $ \Phi$ is holomorphic with respect to a K{\"a}hler
structure if and only if there exists an m-dimensional isotropic subspace $W
\subset {\Bbb C}^{2m} $ such that $\nabla\Phi(x)\in W $ for any $ x \in U $.
\end{proposition}

We will prove that $\Phi = \phi\circ Q_{r}$ is not holomorphic with respect to
any K{\"a}hler structure by showing that there exists no such isotropic
subspace for $\Phi$. It is enough to check that there exist 8 points $
x_{1},\ldots,x_{8} \in  {\Bbb R}^{16} $ such that $ \left\{
\nabla\Phi(x_{1}),\ldots, \nabla\Phi(x_{8})\right \} $ is a linearly
independent set of vectors in $ {\Bbb C}^{16} $, and that there exists another
point $x_{9}\in {\Bbb R}^{16} $  such that $ \nabla\Phi(x_{9})\not \in
span\left\{ \nabla\Phi(x_{1}),\ldots, \nabla\Phi(x_{8}) \right\} $.\\
In fact,
\begin{align}
 \Phi &({z}_{1},{z}_{2},{z}_{3},{z}_{4},{w}_{1},{w}_{2},{w}_{3},{w}_{4}) =
\notag\\
      &({z}_{3}{w}_{1} - {\overline{z}}_{4}{w}_{2} + {z}_{1}{w}_{3} - {z}_{2}
{\overline{w}}_{4})({z}_{4}{w}_{1} + {\overline{z}}_{3}{w}_{2} + {z}_{2}
{\overline{w}}_{3} + {z}_{1}{w}_{4})\notag
\end{align}
and ,
\begin{align}
&\nabla\Phi(x) = \notag\\
&({w}_{3}B +{w}_{4}A , i(w_{3}B + w_{4}A), -{\overline{w}}_{4}B +
{\overline{w}}_{3}A , i (-{\overline{w}}_{4}B + {\overline{w}}_{3}A ) ,
\notag\\
&\  w_{1}B+w_{2}A, i( w_{1}B-w_{2}A), -w_{2}B+w_{1}A, i(w_{2}B+w_{1}A),\notag\\
&({z}_{3}B +{z}_{4}A , i(z_{3}B + z_{4}A), -{\overline{z}}_{4}B +
{\overline{z}}_{3}A , i (-{\overline{z}}_{4}B + {\overline{z}}_{3}A ),\notag\\
&  z_{1}B+z_{2}A, i( z_{1}B-z_{2}A), -z_{2}B+z_{1}A, i(z_{2}B+z_{1}A) )\notag
\end{align}
Where $ A = {z}_{3}{w}_{1} - {\overline{z}}_{4}{w}_{2} + {z}_{1}{w}_{3} -
{z}_{2} {\overline{w}}_{4} , B = {z}_{4}{w}_{1} + {\overline{z}}_{3}{w}_{2} +
{z}_{2} {\overline{w}}_{3} + {z}_{1}{w}_{4} .$\\
We calculate
\begin{align}
&\nabla\Phi(0,0,1,0,1,0,0,1) = (1,i,0,0,0,0,1,i,0,0,1,i,0,0,0,0)\notag\\
&\nabla\Phi(0,0,i,0,1,0,0,1) = (i,-1,0,0,0,0,i,-1,0,0,1,i,0,0,0,0)\notag\\
&\nabla\Phi(1,0,0,0,1,0,1,0) = (0,0,1,i,0,0,1,i,0,0,0,0,0,0,1,i)\notag\\
&\nabla\Phi(i,0,0,0,1,0,1,0) = (0,0,i,-1,0,0,i,-1,0,0,0,0,0,0,-1,-i)\notag\\
&\nabla\Phi(1,0,0,1,1,0,0,0) = (0,0,0,0,1,i,0,0,0,0,-1,-i,1,i,0,0)\notag\\
&\nabla\Phi(1,0,0,1,i,0,0,0) = (0,0,0,0,-1,-i,0,0,0,0,-i,1,i,-1,0,0)\notag\\
&\nabla\Phi(1,0,1,0,1,0,0,0) = (0,0,0,0,0,0,1,i,0,0,1,i,0,0,1,i)\notag\\
&\nabla\Phi(1,0,1,0,i,0,0,0) = (0,0,0,0,0,0,-1,i,0,0,i,-1,0,0,i,-1)\notag\
\end{align}
It is easy to see that these eight vectors are mutually orthogonal in ${\Bbb
C}^{16}$ with respect to $\langle z,w \rangle ^{C} = z_{1}w_{1} + \ldots +
z_{16}w_{16} $ . Hence they are linearly independent.
It is easy to see that $ \nabla\Phi(0,0,1-i,0,1,1,0,0) =
(0,0,0,0,2,-2,-2i,2i,2,2i,2,2i,0,0,0,0)$ does not belong to the subspace
spanned by the above eight vectors.
\end{example}
\begin{remark}
It is known  that the orthogonal multiplications, as maps between Euclidean
spaces, can be used to construct some harmonic maps and morphisms. However the
complete lift of an orthogonal multiplication need not be an orthogonal
multiplication in general as one can check this for the quaternion product.
\end{remark}

{\bf Acknowledgments}: This work was done while the author was visiting
Department of Pure Mathematics, University of Leeds, U.K.. The author is
grateful to the department for the hospitality and generosity, and to E.Loubeau
and many colleaques there for their help and friendship. The author would like
to thank S. Gudmundsson for his generosity in sending the author many useful
references and for his comments which led to an improvement of the original
manuscript. Last and most importantly, the author is greatly indebted to John
C. Wood for making the author's visit to Leeds possible and for guiding the
author in the study of harmonic maps and morphisms. During the preparation of
this work the author benefitted much from his extremely useful comments and
suggestions.

\end{document}